\documentclass[a4paper, 10pt]{article} 

\usepackage{graphicx}
\usepackage{amsmath,amssymb,amsfonts,amsthm,latexsym}
\usepackage{mathrsfs}
\usepackage{textcomp}
\usepackage{tikz-cd}
\usepackage[all]{xy}

\newtheorem{thm}{Theorem}
\newtheorem{cor}[thm]{Corollary}
\newtheorem{lem}[thm]{Lemma}
\newtheorem{prop}[thm]{Proposition}
\theoremstyle{definition}
\newtheorem{defn}[thm]{Definition}
\theoremstyle{remark}

\theoremstyle{remark}
\newtheorem{ex}[thm]{Example}
\theoremstyle{remark}

\theoremstyle{remark}

\newcommand{\F}{\mbox{$\mathcal F$}}

\newcommand{\B}{\mbox{$\mathcal B$}}
\newcommand{\G}{\mbox{$\mathcal{G}$}}

\newcommand\set[1]{\mbox{$\{#1\}$}}
\newcommand\cat[1]{\mbox{$\mathbf{#1}$}}

\newcommand\map[3]{\mbox{$#1:#2\rightarrow #3$}}
\newcommand\inv[1]{\mbox{$#1^{-1}$}}
\newcommand\trip[1]{\mbox{$\mathbb{#1}$}}
\newcommand{\ult}{\mbox{$\mathfrak{U}$}}
\newcommand{\cau}{\mbox{$\mathfrak{C}$}}
\newcommand{\qu}{\mbox{$\mathcal{U}$}}
\newcommand{\prim}{\mbox{$\mathfrak{P}$}}
\newcommand{\id}{\mbox{$\mathfrak{I}$}}
\newcommand{\sob}{\mbox{$\mathfrak{S}$}}
\newcommand\comp[2]{\mbox{$#1\cdot #2$}}
\newcommand\ocomp[2]{\mbox{$#1\circ #2$}}
\newcommand\mon[3]{\mbox{$(#1,#2,#3)$}}
\newcommand\powcat[2]{\mbox{$#1^{\mathbb{#2}}$}}

\newcommand{\arrowtcupp}[2]{\arrow[bend left=50, ""{name=U, below,inner sep=1}]{#1}\arrow[Rightarrow,from=U,to=MU,"#2"]}
\newcommand{\arrowtclow}[2]{\arrow[bend right=50, ""{name=L,inner sep=1}]{#1}\arrow[Rightarrow,from=LM,to=L]{}[]{#2}} 
\newcommand{\arrowtcmid}[2]{\arrow[""{name=MU,inner sep=1},""{name=LM,below,inner sep=1}]{#1}[pos=.1]{#2}}

\title{\textbf{Separated and prime compactifications}\\ 
} 

\author{\textsc{A. Razafindrakoto\thanks{Email: arazafindrakoto@uwc.ac.za - Department of Mathematics and Applied Mathematics - University of the Western Cape}}} 

\date{} 

\begin{document}

\maketitle 

\begin{abstract}
We discuss conditions under which certain compactifications of topological spaces can be obtained by composing the ultrafilter space monad with suitable reflectors. In particular, we show that these compactifications inherit their categorical properties from the ultrafilter space monad. We further observe that various constructions such as the prime open filter monad defined by H. Simmons, the prime closed filter compactification studied by Bentley and Herrlich, as well as the separated completion monad studied by Salbany fall within the same categorical framework.\\
\end{abstract}


\vspace{20pt} 


\section{Introduction}
The purpose of the present work is to present a categorical framework to the approach considered by Salbany in his article \cite{Sal00} to describe certain compactifications. Salbany showed in \cite{Sal00} that standard compactifications such as the $T_0$ stable compactification (\cite{Smy}) and the \v{Cech}-Stone compactification can be seen as quotients of the space of ultrafilters on the designated spaces. The construction of the space of ultrafilters itself yields a monad, called {\em ultrafilter space monad}, on the category of topological spaces as shown in \cite{Low}. This monad, as we shall show, plays a similar role to that of the completion monad studied by Salbany in the article \cite{Sal82} in the context of quasi-uniform spaces. Under some mild conditions a suitable reflector composed with the ultrafilter space monad gives the desired compactification. In particular, many of the results in \cite{Sal00} can be obtained from this observation.

The fact that reflective subcategories play an important role in Topology is well-known as is visible in the work of Herrlich's \cite{Her}. On the other hand, monads have grown increasingly important in Topology especially through the influential work of Barr (\cite{Bar}), Manes (\cite{Man}), Br\"ummer (\cite{Bru79}), Salbany (\cite{Sal84}), Simmons (\cite{Sim}) and Wyler (\cite{Wyl74}). In Pointfree Topology, the presence of monads can be explicitly observed in Johnstone's book \cite{Joh} and in Banaschewski and Br\"ummer's paper \cite{BanBru}. Recently, monads were considered as fundamental building blocks in the study of lax algebras (\cite{Low}) and contributed to further applications of categorical techniques in Topology. Following these various developments, we show that monads play some significant role in certain familiar constructions. Indeed, by observing that idempotent monads are essentially reflectors, we can consider the algebras of a monad  as sitting in a generalised reflective subcategory. This is especially important for compactifications such as the closed prime filter compactification (\cite{BenHer}) which is not idempotent. 

In this paper, we consider a monad $T$ together with a reflector $R$, and present conditions under which the composition $R\cdot T$ itself becomes a reflector. In categorical parlance, we provide conditions that permit the distributivity of $T$ over $R$, a fact that expresses the harmonious interplay between separation and compactness. In the examples of compactifications furnished here, these conditions are intimately linked to the fact that continuous maps between compact and separated spaces must be proper (or perfect). Among other things, we show that $R\cdot T$ can be seen as the universal reflection of $T$ among suitable monads on the ambient category. Together with the {\em Boolean Ultrafilter Theorem}\footnote{Every proper filter on a given Boolean algebra is contained in an ultrafilter.} (BUT), this is used to deduce that some classical instances of compactifications are suitable quotients of the ultrafilter space.

The paper is organised in a simple manner. After providing the strict necessary background in the second section, we review Salbany's construction of the ultrafilter space as well as Simmons' prime open filter monad in the third section. This section not only motivates the categorical approach but provides technical details that are similar to the proofs that shall be omitted in the example on prime closed filter monad. The main results are discussed in Section 4, followed with some fundamental examples. The examples  are admittedly not exhaustive but are chosen to illustrate the theory. For instance, spectral spaces, bitopological spaces and the Samuel compactification which deserve to be treated, are unfortunately left out. The abstract language is then necessary not only for the duality it provides (e.g. the case for frames) and its capacity to handle similar constructions in different settings, but also for its role as a guide for future examples.

\section{Background and preliminaries}

\subsection{Categorical background}
The composition of morphisms $f$ and $g$ in a category shall be denoted by \comp{g}{f} or just simply by $gf$. Given two pairs of parallel functors \map{F,G}{\cat{C}}{\cat{D}} and \map{H,K}{\cat{D}}{\cat{E}} and two natural transformations $\alpha:F\to G$ and $\beta:H\to K$, the natural transformation which is given by the horizontal composition $\beta\circ\alpha:HF\to KG$ is defined as follows: for each $X$ in \cat{C}
\begin{center}
$(\beta\circ\alpha)_X=\beta_{GX}H(\alpha_X)=K(\alpha_X)\beta_{FX}$.
\end{center}
In most cases, the variable $X$ will be omitted and we shall simply write $\beta\circ\alpha=\beta G\cdot H\alpha=K\alpha\cdot\beta F$.
If \map{L}{\cat{C}}{\cat{D}} is another functor and $\delta:G\to L$ another natural transformation, then the vertical composition $\comp{\delta}{\alpha}:F\to L$ is simply defined as $\comp{\delta_X}{\alpha}_X$ for each $X\in\cat{C}$. Recall that 
\begin{lem}
\label{lem: middle-interchange law}
(Middle-Interchange Law) Given natural transformations
\[
\begin{tikzcd}[column sep=2cm]
 \cat{C}   \arrowtcmid{r}{} \arrowtcupp{r}{\alpha}\arrowtclow{r}{\alpha\smash'} & \cat{D} \arrowtcmid{r}{} \arrowtcupp{r}{\beta}\arrowtclow{r}{\beta\smash'} & \cat{E}
\end{tikzcd}\]
we have $(\beta\smash'\circ\alpha\smash')\cdot(\beta\circ\alpha)=(\beta\smash'\cdot\beta)\circ(\alpha\smash'\cdot\alpha)$. (See \cite{Low,Mac}.)
\end{lem}

A subcategory \cat{D} of \cat{C} is said to be {\em reflective} if it is full and the inclusion functor $\cat{D}\to\cat{C}$ admits a left adjoint \map{R}{\cat{C}}{\cat{D}}. If \map{r}{1}{R} is the unit of $R$, then the pair $(R,r)$ is called a {\em reflector}. We say that the pair $(R,r)$ is an {\em epireflector} if each \map{r_X}{X}{RX} is an epimorphism in the category \cat{C}. We usually write $R(\cat{C})$ for \cat{D} and just consider \map{R}{\cat{C}}{\cat{C}} as an endofunctor.

A {\em monad} $\mathbb{T}$ on a category \cat{C} is a triple \mon{T}{\mu}{\eta}, where \map{\mu}{TT}{T} and \map{\eta}{1}{T} are natural transformations satisfying the identities 
\begin{center}
$\comp{\mu}{T\mu}=\comp{\mu}{\mu T}$ and $\comp{\mu}{\eta T}=\comp{\mu}{T\eta}=1_T$.
\end{center}
A morphism between two monads $\trip{T}=\mon{T}{\mu}{\eta}$ and $\trip{M}=\mon{M}{m}{e}$ is a natural transformation \map{\alpha}{T}{M} satisfying $\comp{\alpha}{\eta}=e$ and $\comp{\alpha}{\mu}=\comp{m}({\ocomp{\alpha}{\alpha}})$. A \trip{T}-algebra (or an {\em Eilenberg-Moore algebra}) is a pair $(X,a)$, where $X\in\cat{C}$ and \map{a}{TX}{X} a morphism called {\em structure morphism} such that 
\begin{center}
$\comp{a}{Ta}=\comp{a}{\mu_X}$ and $\comp{a}{\eta_X}=1_X$.
\end{center}
Note that $(TX,\mu_X)$ is the free \trip{T}-algebra over $X$. If $(X,a)$ and $(Y,b)$ are \trip{T}-algebras, then a \trip{T}-algebra homomorphism \map{f}{(X,a)}{(Y,b)} is a morphism \map{f}{X}{Y} in \cat{C} such that $\comp{f}{a}=b\cdot Tf$. We note that structure morphisms \map{a}{TX}{X} are \trip{T}-algebra homomorphisms. The category of \trip{T}-algebras and \trip{T}-algebra homomorphisms are denoted by \powcat{\cat{C}}{T}. The forgetful functor $\map{\powcat{G}{T}}{\powcat{\cat{C}}{T}}{\cat{C}}:(X,a)\mapsto X$ admits a left adjoint $\map{\powcat{F}{T}}{\cat{C}}{\powcat{\cat{C}}{T}}:X\mapsto (TX,\mu_X), f\mapsto Tf$. The unit of this adjunction is given by $\eta_X:X\to \powcat{G}{T}\powcat{F}{T}X$ and the co-unit is provided by structure maps $\varepsilon_{(X,a)}:\powcat{F}{T}\powcat{G}{T}(X,a)\to(X,a)$. We wish to emphasize the following observations.

\begin{lem}
\label{lem: uniqueness of structure map}
(See \cite[Section 4, Item 12]{Man}.) For each morphism \map{f}{X}{\powcat{G}{T}(Y,b)} in \cat{C}, there is a unique \trip{T}-algebra homomorphism \map{\bar{f}}{(TX,\mu_X)}{(Y,b)}, namely $\bar{f}=\comp{b}{Tf}$, such that $\powcat{G}{T}\bar{f}\cdot\eta_X=f$. In particular there is at most a unique \trip{T}-algebra homomorphism \map{a}{TX}{X} such that $a\eta_X=1_X$.
\end{lem}

\begin{lem}
\label{lem: equivalence reflective subcategories and idempotent monads}
(See \cite[Corollary 4.2.4]{Bor2}.) Reflective subcategories of \cat{C} coincide, up to equivalences of categories, with categories \powcat{\cat{C}}{T} of \trip{T}-algebras for the idempotent monads \trip{T} on \cat{C}.
\end{lem}

\begin{ex}
\label{ex: ultrafilter monad on sets}
Consider the category \cat{Set} of sets and functions, and for each set $X$ the collection $\ult X$ of ultrafilters on $X$. To each function \map{f}{X}{Y} is assigned a function $\ult f$ defined by
\begin{center}
$\ult f(\F)=\{B\subseteq Y\ |\ \inv{f}(B)\in\ult X\}$.
\end{center}
Thus \ult\ is an endofunctor and it forms a monad $\mathbb{U}=(\ult,\mu,\eta)$ on \cat{Set} with the multiplication $\mu$ and the unit $\eta$ defined by:
\begin{center}
$\mu_X(\mathfrak{X})=\{A\subseteq X\ |\ A^*\in\mathfrak{X}\}$ and $\eta_X(x)=\{A\subseteq X\ |\ x\in A\}$
\end{center}
where $A^*=\{\F\in\ult X\ |\ A\in\F\}$,\ $\mathfrak{X}\in\ult\ult X$ and $x\in X$. It was shown by E. Manes \cite{Man} that \powcat{\cat{Set}}{U} is isomorphic to the category of compact Hausdorff spaces and continuous maps.
\end{ex}

\subsection{Stably compact spaces}

A given topological space $X$ comes with a natural pre-ordering through the {\em specialisation order}: $x\leq y$ if and only if $x\in\overline{\{y\}}$ or $\overline{\{x\}}\subseteq\overline{\{y\}}$.  $T_0$ spaces are those spaces for which the equality $\overline{\{x\}}=\overline{\{y\}}$ implies $x=y$. The specialisation order is then anti-symmetric in a $T_0$ space.

A nonempty closed set $G$ is called irreducible if the inclusion $G\subseteq F_1\cup F_2$ implies $G\subseteq F_1$ or $G\subseteq F_2$ whenever $F_1$ and $F_2$ are closed. A space $X$ is {\em sober} if all irreducible closed sets in $X$ are of the form $\overline{\{x\}}$, where such $x\in X$ must be unique. Thus sobriety implies the $T_0$ separation axiom. If this property of uniqueness is not granted, then $X$ shall be called {\em weakly sober} (\cite{Low}).

For two open sets $O$ and $U$ in a topological space $X$, we say that $O$ is {\em relatively compact} in $U$ - and we write $O\ll U$, when every open covering of $U$ admits a finite sub-covering of $O$. A space $X$ is said to be {\em stable} if the order relation $\ll$ is finitely multiplicative, i.e. $X\ll X$ and for any open sets $O, U, W$ and $V$, if $O\ll U$ and $W\ll V$ then $O\cap W\ll U\cap V$. Finally, we say that $X$ is {\em locally compact} if compact neighbourhoods form a base for the neighbourhood system of the space $X$.

\begin{defn}
A space $X$ that is locally compact, stable and weakly sober shall be called a {\em Salbany space}. $X$ is called {\em stably compact} if $X$ is a $T_0$ Salbany space.
\end{defn}
Salbany spaces are called {\em quasi-stably compact spaces} in \cite{Bez} and Salbany called them {\em stably compact spaces} in \cite{Sal00}. These are precisely the {\em representable topological spaces} in \cite{Low}. On the other hand, stably compact spaces went under various names: {\em strongly sober locally compact} (in \cite{KueBru}) and {\em stably quasi-compact} (in \cite{Hof84}). Note that these nomenclatures may be potentially confusing over time: for instance {\em stably quasi-compact} and {\em quasi-stably compact spaces} may be perceived as the same class of spaces. We find it convenient to adopt a new name that does not include the prefixes `quasi' and `stably'.

A subset $A\subseteq X$ such that $A=\bigcap\{O\subseteq X\ |\ O\text{ open and }A\subseteq O\}$ is called {\em saturated}. A continuous map \map{f}{X}{Y} between two Salbany spaces is said to be {\em proper}\footnote{Also called {\em perfect}, these are the pseudo-homomorphisms in \cite{Low}.} if for any compact  saturated set $K\subseteq Y$, the inverse image $\inv{f}(K)$ is compact in $Y$. The category of Salbany spaces and proper maps is denoted by \cat{Sal} and that of stably compact spaces and proper maps by \cat{Stb}.

\begin{ex}
\label{ex: salbany space that is not stably compact}
Salbany spaces which are not stably compact can be obtained by ``inserting'' points into the irreducible closed subsets. An example in the finite case is the set $X=\set{1,2,3}$ with the topology $\tau=\set{\emptyset,X,\set{1}}$. $(X,\tau)$ is easily seen to be stable, locally compact and weakly sober with irreducible closed subsets \set{X,\set{2,3}}.
\end{ex}

\paragraph*{Patch topology.}
If $(X,\tau)$ is a stably compact space, then one can form a new topology $\tau^c$ with the collection $\{K\subseteq X\ |\ K\text{ compact saturated}\}$ as basic closed sets. The join topology $\pi=\tau\vee\tau^c$ is called the {\em patch topology} on $(X,\tau)$. In principle, the patch topology can be constructed for any given topological space. However, in the case where $X$ is stably compact, the topology $\pi$ is compact and Hausdorff (\cite{Bez}). It can be shown that if $\sigma$ is another compact Hausdroff topology on $X$ such that $\tau\subseteq\sigma$, then $\pi\subseteq\sigma$. Consequently we have that
\begin{lem}
\label{lem: CHaus coreflective in Stb}
(\cite{Esc}) Compact Hausdorff spaces are coreflective in the category \cat{Stb}.
\end{lem} 
It is known that \cat{Stb} is equivalent to the category of {\em ordered compact Hausdorff spaces} and monotone continuous maps. These spaces are discussed by various authors in \cite{Bez,Jun,Low,Hof84,Wyl84}.

\section{The ultrafilter space monad}
\label{The ultrafilter space monad}

The ultrafilter monad $\mathbb{U}=(\ult,\mu,\eta)$ on \cat{Set} from Example \ref{ex: ultrafilter monad on sets} can be lifted to \cat{Top} by defining the topology with basis\footnote{Various equivalent constructions are given in \cite{Bez} and \cite{Low}.}
\begin{center}
$\{O^*\ |\ O\subseteq X\text{ open}\}$ where $O^*=\{\F\in\ult X\ |\ O\in\F\}$
\end{center}
on the set of ultafilters $\ult X$ for each space $X$. That $\mu$ and $\eta$ are continuous follow in a straightforward way from the facts that $\inv{\mu_X}(O^*)=(O^*)^*$ and $\inv{\eta_X}(O^*)=O$ for each open set $O\subseteq X$. This lifting into the category \cat{Top} shall still be denoted with the same symbols $\mathbb{U}=(\ult,\mu,\eta)$ since we will not concern ourselves with \cat{Set}, and we will refer to $\mathbb{U}$ as the {\em ultrafilter space monad}. The following were shown by Salbany in \cite{Sal00}:

\begin{prop}
\label{prop: existence of continuous retraction}
(See \cite[Theorem 1, Propositions 2, 3 and 5]{Sal00}.) A space $X$ is a Salbany space if and only if the embedding $\eta$ admits a continuous retraction \map{r}{\ult X}{X}. The space $\ult X$ is a salbany space and the embedding \map{\eta}{X}{\ult X} is dense in the patch topology of $\ult X$. 
\end{prop}

Note that a given continuous retraction as in Proposition \ref{prop: existence of continuous retraction} is not necessarily an $\mathbb{U}$-homomorphism since it is not unique as already pointed out in \cite[Proposition 3]{Sal00}. Although Salbany showed that any retraction $r:\ult X\to X$ is proper (\cite[Proposition 4]{Sal00}), this is not enough to conclude that \cat{Sal} is the category of $\mathbb{U}$-algebras. Instead, Salbany spaces and proper maps form the category of {\em pseudo-algebras and pseudo-homomorphisms} as shown in \cite[Theorem 5.7.2]{Low}. Algebras of the ultrafilter space monad were identified up to equivalence in \cite[Lemma 5.6.1]{Low} and fully  characterised in \cite{Bez}.

\begin{thm}
\label{thm: algebras of the ultrafilter space monad}
(\cite[Lemma 4.12 and Theorem 4.15]{Bez}) The algebras \powcat{\cat{Top}}{U} of $\mathbb{U}$ are exactly the bitopological spaces $(X,\tau,\pi)$ together with maps \map{f}{(X,\tau,\pi)}{(Y,\tau',\pi')} such that both \map{f}{(X,\tau)}{(X,\tau')} and \map{f}{(X,\pi)}{(X,\pi')} are continuous and such that:
\begin{enumerate}
\item $(X,\tau)$ is a Salbany space and $(X,\pi)$ is compact Hausdorff;
\item $\tau\subseteq\pi$ and every compact saturated subset of $(X,\tau)$ is closed in $(X,\pi)$.
\end{enumerate}
\end{thm}
The topology $\pi$ in Theorem \ref{thm: algebras of the ultrafilter space monad} owes its existence and uniqueness to the fact that the underlying set-monad of $\mathbb{U}$ admits precisely the compact Hausdorff spaces as algebras (\cite{Man}). It is important to point out in this case that if $\tau$ is not sober, then $\pi$ does not necessarily agree with the patch topology construction $\tau\vee\tau^c$. (See Example \ref{ex: salbany space that is not stably compact} for instance.) The category \powcat{\cat{Top}}{U} is equivalent to the category of {\em pre-ordered compact Hausdorff spaces} and monotone continuous maps. This equivalence is discussed in \cite[Theorem 5.5]{Bez} and \cite[Lemma 5.6.1]{Low}.

Our task in the remaining part of this section is to describe the passage from \powcat{\cat{Top}}{U} to \cat{Stb} through a reflector. One of the two monads described by Simmons in \cite{Sim} is the {\em prime open filter} monad $\mathbb{S}=(\Sigma,m,e)$, where $\Sigma X=\set{\G\text{ filter } |\ \G\subseteq\tau\text{ is prime with respect to open subsets}}$ for each space $X$, with
\begin{center}
$m_X(\mathfrak{X})=\set{O\text{ open }\ |\ O^{\times}\in\mathfrak{X}}$ and $e_X(x)=\{O\ |\ O\text{ open and }x\in O\}$,
\end{center}
and where $O^{\times}=\{\G\in\Sigma X\ |\ O\in\G\}$, $x\in X$ and $\mathfrak{X}\in\Sigma\Sigma X$. The topology on $\Sigma X$ has as a basis the sets \set{O^{\times}\ |\ O\subseteq X\text{ open}}. For each continuous map $f$, $\Sigma f$ acts as a restriction of $\ult f$. It was shown by Simmons in \cite[Theorem 1.3 and Lemma 3.12]{Sim} that $\cat{Top}^{\trip{S}}\cong\cat{Stb}$. 
\begin{prop}
\label{prop: monad morphism from prime to open prime}
For each space $(X,\tau)$, the map \map{\alpha_X}{\ult X}{\Sigma X} that takes each ultrafilter \F\ to the prime open filter $\F\cap\tau$ is part of a morphism of monads \map{\alpha}{\ult}{\Sigma}.
\end{prop}

\begin{proof}
Each $\alpha_X$ is clearly well-defined. If $O^{\times}$ is a basic open set in $\Sigma X$, then 
\begin{center}
$\inv{\alpha_X}(O^{\times})=\set{\F\in\ult X\ |\ \F\cap\tau\in O^{\times}}=\set{\F\in\ult X\ |\ O\in\F\cap\tau}=O^*$.
\end{center}
If \map{f}{X}{Y} is a continuous map, then we easily have $\Sigma f\cdot\alpha_X=\alpha_Y\cdot\ult f$. That $e=\alpha\cdot\eta$ is clear. It remains to show that $m\cdot(\alpha\circ\alpha)=\alpha\cdot\mu$. First note that we have $\alpha\circ\alpha=\alpha\Sigma\cdot\ult\alpha$ so that for each $\mathfrak{X}\in\ult\ult X$
\begin{align*}
(\alpha\Sigma\cdot\ult\alpha)(\mathfrak{X}) &= \alpha\Sigma\left(\set{\chi\subseteq\Sigma X\ |\ \inv{\alpha}(\chi)\in\mathfrak{X}}\right)\\
&= \set{\chi\subseteq\Sigma X\ |\ \inv{\alpha}(\chi)\in\mathfrak{X}}\cap\set{\chi\subseteq\Sigma X\ |\ \chi\text{ open}}\\
&=\set{\chi\subseteq\Sigma X\ |\ \inv{\alpha}(\chi)\in\mathfrak{X}\text{ and }\chi\text{ open}}.
\end{align*}
Now $(m\cdot(\alpha\circ\alpha))(\mathfrak{X})=(m\cdot\alpha\Sigma\cdot\ult\alpha)(\mathfrak{X})=\set{O\text{ open}\ |\ \inv{\alpha}(O^{\times})\in\mathfrak{X}}$. On the other hand
\begin{center}
$(\alpha\cdot\mu)(\mathfrak{X})= \set{W\subseteq X\ |\ W^*\in\mathfrak{X}}\cap\tau= \set{O\text{ open}\ |\ \inv{\alpha}(O^{\times})\in\mathfrak{X}}.$
\end{center}
We then have equality.
\end{proof}
We now show that $\Sigma X$ appears as the $T_0$ reflection of $\ult X$.

\begin{prop}
\label{prop: reflection of prime monad is open prime monad}
Let \map{(R,r)}{\cat{Top}}{\cat{Top}} be the $T_0$ reflector. For each space $X$, there is an isomorphism $r_X\cong\alpha_X$ in a sense that there is an homeomorphism \map{\varphi}{R\ult X}{\Sigma X} with $\varphi\cdot r_X=\alpha_X$.
\end{prop}

\begin{proof}
The existence of $\varphi$ as a continuous function follows from the fact that $\Sigma X$ is a $T_0$ space (\cite{Sim}). For each equivalence class $r(\F)\in R\ult X$ we have $\varphi(r(\F))=\F\cap\tau$. Since $r(\F)=r(\G)$ if and only if $\overline{\set{\F}}=\overline{\set{\G}}$, and that $\F\in O^*$ if and only if $O\in\F$, $\varphi$ is injective. That $\varphi$ is surjective follows from BUT. Finally as $\alpha$ is initial and $r$ is surjective, $\varphi$ is initial and hence open.
\end{proof}
Proposition \ref{prop: monad morphism from prime to open prime} and Proposition \ref{prop: reflection of prime monad is open prime monad} imply that \map{r\ult}{\ult}{R\ult} is a morphism of monads. In fact this is a consequence of a much deeper aspect of $(R,r)$ and $\mathbb{U}$, namely that the composition $R\cdot\ult$ is a monad composition $\mathbb{R}\circ\mathbb{U}$.

\begin{prop}
\label{prop: t0 quotient maps are proper}
Let $(R,r)$ be the $T_0$ reflector on \cat{Top}. For each space $X$, $r_X$ is a proper map.
\end{prop}

\begin{proof}
For arbitrary open sets $O\subseteq X$ and $U\subseteq RX$,\ $\inv{r}(r(O))=O$ and $r(\inv{r}(U))=U$. This establishes a lattice isomorphism between the topology of $X$ and that of $RX$. Thus $K\subseteq RX$ compact is equivalent to $\inv{r}(K)\subseteq X$ being compact.
\end{proof}
Now if \map{f}{X}{Y} is a proper map where $Y$ is a $T_0$ space, then the unique continuous map \map{\varphi}{RX}{Y} such that $\varphi\cdot r=f$ is proper. Therefore
\begin{cor}
\cat{Stb} is reflective in both \cat{Sal} and \powcat{\cat{Top}}{U}.
\end{cor}
If \map{I}{\cat{Stb}}{\powcat{\cat{Top}}{U}} is the full inclusion, then adjunction between \cat{Top} and \cat{Stb} is decomposed as follows\footnote{This adjunction is described as a reflection by Simmons in \cite[Theorem 1.3]{Sim}. This notion of reflectivity does not include idempotency of the reflector. In our case, since the subcategory \cat{Stb} is not full, $\Sigma$ cannot be idempotent.}
$$\xymatrix{
\cat{Top} \ar@<1.2ex>[rr]^{\powcat{F}{U}}_{\perp} & & \powcat{\cat{Top}}{U} \ar@<1.5ex>[ll]^{\powcat{G}{U}} \ar@<1.1ex>[rr]^R
& & \cat{Stb} \ar@<1ex>[ll]^I_{\perp} 
}$$
Salbany has shown in \cite{Sal00} - citing \cite{Smy} and \cite{Sim}, that if $X$ is $T_0$ space, then \map{r\ult\cdot\eta}{X}{R\ult X} or equivalently \map{e}{X}{\Sigma X} is the $T_0$ stable compactification of $X$.

\section{Distributivity of a monad over an epireflector}

We fix a reflector $\mathbb{R}=(R,r)$ and a monad $\mathbb{T}=(T,\mu,\eta)$ on a given category \cat{C}. Let us start with the following result which is from \cite[Lemma II.3.8.1]{Low} and \cite[Propositions II.3.8.2 and II.3.8.4]{Low}. (See also \cite[Theorem 2.4.2]{ManMul}.)

\begin{thm}
\label{thm: lifting of a monad}
(See \cite{ManMul} and \cite{Low}.) The composite functor $RT$ underlies a composite monad $\mathbb{R\circ T}=(RT,w,r\circ\eta)$ if and only if there is a lifting of $R$ on \powcat{\cat{C}}{T}, i.e. there is a monad $\mathbb{\tilde{R}}=(\tilde{R},n,d)$  such that the diagram
$$\xymatrix{
\powcat{\cat{C}}{T} \ar[r]^{\tilde{R}} \ar[d]_{\powcat{G}{T}} & \powcat{\cat{C}}{T} \ar[d]^{\powcat{G}{T}}\\
\cat{C} \ar[r]_R & \cat{C}
}$$
commutes, $\powcat{G}{T}n=1_R\powcat{G}{T}$ and $\powcat{G}{T}d=r\powcat{G}{T}$. Under these conditions, \map{rT}{T}{RT} is a monad morphism and $(\powcat{\cat{C}}{T})^{\mathbb{\tilde{R}}}\cong\powcat{\cat{C}}{R\circ T}$. 
\end{thm}
This means that for each $\mathbb{T}$-algebra $(X,a)$, one has $\tilde{R}(X,a)=(RX,\tilde{R}a)$ for some structure morphism \map{\tilde{R}a}{TRX}{RX}. The multiplication $w$ is obtained  as the reflection $R\tilde{R}\mu$. Since {\powcat{G}{T}} forgets algebra structures, the underlying \cat{C}-morphisms of $n$ and $d$ coincide with $1_R$ and $r$ respectively. Thus also $n=1_{\tilde{R}}$, and since \map{d}{(X,a)}{\tilde{R}(X,a)} is a morphism in \powcat{\cat{C}}{T}, the necessary condition above can be exclusively expressed through $\mathbb{T}$ and $\mathbb{R}$.

\begin{lem}
\label{lem: lifting a reflector as a monad}
The composite functor $RT$ is part of a monad $\mathbb{R\circ T}$ if for each $\mathbb{T}$-algebra $(X,a)$, the morphism \map{r}{X}{RX} is a $\mathbb{T}$-algebra homomorphism, i.e. there is a morphism \map{b}{TRX}{RX} making $(RX,b)$ into a \trip{T}-algebra and such that $b\cdot Tr=r\cdot a$.
\end{lem}

\begin{prop}
\label{prop: algebras of composite are reflective}
$(\powcat{\cat{C}}{T})^{\mathbb{\tilde{R}}}$ is isomorphic to a reflective subcategory of \powcat{\cat{C}}{T}.
\end{prop}

\begin{proof}
If $\mathbb{R\circ T}$ is a monad, then the lifting \map{\tilde{R}}{\powcat{\cat{C}}{T}}{\powcat{\cat{C}}{T}} produces an adjunction
$$\xymatrix{
\powcat{\cat{C}}{T} \ar@<1ex>[rr]^{\powcat{F}{\tilde{R}}}_{\perp} & &  (\powcat{\cat{C}}{T})^{\mathbb{\tilde{R}}}\ar@<1.3ex>[ll]^{\powcat{G}{\tilde{R}}} 
}$$
Since $R$ is idempotent, we have $\tilde{R}(\tilde{R}(X,a))\cong \tilde{R}(X,a)$ for any $\mathbb{T}$-algebra $(X,a)$. Thus $\tilde{R}$ is also idempotent.
\end{proof}

\begin{cor}
\label{cor: equivalence reflectivity and algebraicity or reflection}
Suppose that $(R,r)$ is a reflector. $\powcat{\cat{C}}{R\circ T}$ is isomorphic to a reflective subcategory of \powcat{\cat{C}}{T} if and only if  \map{r}{X}{RX} is a \trip{T}-algebra homomorphism for each \trip{T}-algebra $(X,a)$.
\end{cor}

A potential difficulty associated to Theorem \ref{thm: lifting of a monad} and Lemma \ref{lem: lifting a reflector as a monad} is that the algebras \powcat{\cat{C}}{T} ought to be identified. In certain categories, the properties of the free algebras $(TX,\mu_X)$ can be understood before the algebras are characterised, and this is certainly the case for Salbany's paper \cite{Sal00}. In what follows, we shall provide conditions on the free algebras and on $R$ that will render the lifting in Theorem \ref{thm: lifting of a monad} possible without fully understanding the algebras of $\mathbb{T}$.

\begin{lem}
\label{lem: transition epimorphism}
If $(R,r)$ is an epireflector, then the following are equivalent:
\begin{enumerate}
\item For each parallel pair of morphisms \map{f,g}{TX}{RY}, the equations $f\eta=g\eta$ and $f=g$ are equivalent.
\item \map{R\eta}{RX}{RTX} is an epimorphism in $R(\cat{C})$.
\end{enumerate}
\end{lem}
\begin{lem}
\label{lem: preservation of algebras}
Suppose that $(R,r)$ is an epireflector and that $R\eta$ is an epimorphism in $R(\cat{C})$. If for a $\mathbb{T}$-algebra $(X,a)$,\ $RX$ admits a \trip{T}-structure morphism $b$, then \map{r}{X}{RX} is a $\mathbb{T}$-algebra homomorphism.
\end{lem} 

\begin{proof}
Since the outer diagram of
$$\xymatrix{
X\ar[r]^{\eta} \ar[d]_r & TX\ar[r]^a \ar[d]^{Tr} & X\ar[d]^r \\
RX\ar[r]_{\eta R} & TRX\ar[r]_b & RX
}$$
commutes, we have $(ra)\eta=(bTr)\eta$ which is equivalent to $ra=bTr$.
\end{proof}
Note that with the assumptions of Lemma \ref{lem: preservation of algebras}, there is at most one \cat{C}-morphism \map{b}{TRX}{RX} such that $b\cdot\eta R=1_R$ and when this morphism exists, then $(RX,b)$ is a \trip{T}-algebra. In fact it is enough for $\eta RT$ to split for $RX$ to become a \trip{T}-algebra, as shown below.

\begin{thm}
\label{thm: preservation of free algebras is enough}
Suppose that $(R,r)$ is an epireflector and $R\eta$ an epimorphism in $R(\cat{C})$. If for each $X\in\cat{C}$,\ $\eta_{RTX}$ is a split monomorphism, then for each $\mathbb{T}$-algebra $(X,a)$,\ $RX$ admits a \trip{T}-structure morphism $b$.
\end{thm}

\begin{proof}
Let $\beta$ be a morphism such that $\beta\cdot\eta{RT}=1_{RT}$ and define $b:=Ra\cdot \beta\cdot TR\eta$. Note first that $TR\eta\cdot\eta R=\eta RT\cdot R\eta$ since $\eta$ is a natural transformation. We have $b\cdot\eta R=(Ra\cdot \beta\cdot TR\eta)\cdot \eta R=Ra\cdot \beta\cdot \eta RT\cdot R\eta=Ra\cdot R\eta=R(a\cdot\eta)=1_R$. This also makes the outer diagram of
$$\xymatrix{
TRX\ar[r]^{\eta TR} \ar[d]_b & TTRX\ar[r]^{\mu R} \ar[d]^{Tb} & TRX\ar[d]^b \\
RX\ar[r]_{\eta R} & TRX\ar[r]_b & RX
}$$
commutative and $b$ a $\mathbb{T}$-algebra homomorphism.
\end{proof}

\begin{cor}
\label{cor: strong lifting of a reflector as a monad}
Suppose that $(R,r)$ is an epireflector and $R\eta$ an epimorphism in $R(\cat{C})$. If for each $X\in\cat{C}$,\ $\eta_{RTX}$ is a split monomorphism, then the composite functor $RT$ forms a monad $\mathbb{R\circ T}$.
\end{cor}

We now discuss certain aspects of the behaviour of \powcat{\cat{C}}{R\circ T} with respect to \powcat{\cat{C}}{T} and $R(\cat{C})$. The following Lemma is a consequence of Theorem \ref{thm: lifting of a monad}.

\begin{lem}
\label{lem: an rt-algebra is a t-algebra}
If $X$ is an \trip{R\circ T}-algebra, then it is a \trip{T}-algebra.
\end{lem}

\begin{proof}
If \map{\alpha}{RTX}{X} is an \trip{R\circ T}-structure morphism, then \comp{\alpha}{rT}\ is a  \trip{T}-structure morphism.
\end{proof}

\begin{lem}
\label{lem: properties related to idempotency}
Suppose that $(R,r)$ is an epireflector.
\begin{enumerate}
\item If $X$ is an \trip{R\circ T}-algebra, then $X\cong RX$.
\item If $R\eta$ is an epimorphism in $R(\cat{C})$, then $RT$ is idempotent.
\end{enumerate}
\end{lem}

\begin{proof}
For the first statement, if \map{\alpha}{RTX}{X} is an \trip{R\circ T}-structure morphism, then we have $(\comp{\alpha}{R\eta})\cdot r=\comp{\alpha}{(rT\cdot\eta)}=1_X$. Therefore $r$, being an epimorphism, is an isomorphism. For the second statement, it is enough to show that if $X$ is \trip{R\circ T}-algebra, then $X\cong RTX$. If \map{\alpha}{RTX}{X} is an \trip{R\circ T}-structure morphism, then
\begin{center}
$(\ocomp{r}{\eta})\cdot\alpha\cdot rT\cdot\eta=\comp{\comp{(\ocomp{r}{\eta})}{\alpha}}{(\ocomp{r}{\eta})}=\ocomp{r}{\eta}=rT\cdot\eta$.
\end{center}
Since $rT$ is an epimorphism and $R\eta$ is an epimorphism in $R(\cat{C})$, we have $(\ocomp{r}{\eta})\cdot\alpha=1_{RTX}$.
\end{proof}

\begin{prop}
\label{prop: properties of the rt-algebras as reflective subcategories}
If \trip{R\circ T} is a monad, $(R,r)$ an epireflector and $R\eta$ an epimorphism in $R(\cat{C})$, then  \powcat{\cat{C}}{R\circ T} is isomorphic to a reflective subcategory of $R(\cat{C})$ and \cat{C}. 

\end{prop}

\begin{proof}
Clearly if $RT$ is idempotent then \powcat{\cat{C}}{R\circ T} is a reflective subcategory of \cat{C} (Lemma \ref{lem: equivalence reflective subcategories and idempotent monads}). Now, if $RY$ is an \trip{R\circ T}-algebra (as in Lemma \ref{lem: properties related to idempotency}.1) and \map{f}{RX}{RY} a morphism in $R(\cat{C})$, then there is a unique \trip{R\circ T}-algebra homomorphism \map{\varphi}{RTRX}{RY} such that $\comp{\varphi}{(\ocomp{r}{\eta})R}=f$. Since the inclusion $R(\cat{C})\to\cat{C}$ is full, both $\varphi$ and $(\ocomp{r}{\eta})R$ are in $R(\cat{C})$.
\end{proof}

A test for (non-)idempotency has been given by Fakir in \cite{Fak} and it allows for a partial converse of Proposition \ref{prop: properties of the rt-algebras as reflective subcategories}.

\begin{prop}
\label{prop: Fakir test of idempotency}
(See \cite[Proposition 1]{Fak}.) Suppose that $\trip{W}=\mon{\Sigma}{m}{e}$ is an idempotent monad on \cat{C}. If \map{e_X}{X}{\Sigma X} is monomorphism for each object $X$, then it is an epimorphism.
\end{prop}

\begin{cor}
\label{cor: Fakir test of idempotency}
Suppose that \powcat{\cat{C}}{R\circ T} is reflective in $R(\cat{C})$ with reflector $(RT,\ocomp{r}{\eta})$. If for any object $X$, the unit \map{(\ocomp{\eta}{r})R}{RX}{RTRX} is a monomorphism, then it is an epimorphism.
\end{cor}

Now, consider another monad $\trip{M}=\mon{M}{n}{e}$ on \cat{C}.

\begin{thm}
\label{thm: universality of the separation}
Assume that $(R,r)$ is an epireflector and that $T$ preserves epimorphisms. Suppose that \trip{R\circ T} is a monad and that for each object $X\in\cat{C}$,\ \ $MX\in R(\cat{C})$. Then for any morphism of monads \map{\gamma}{T}{M}, there is a unique morphism of monads \map{\lambda}{\trip{R\circ T}}{\trip{M}} such that $\lambda\cdot rT=\gamma$.
\end{thm}

\begin{proof}
For each object $X$, there is a unique morphism \map{\lambda_X}{RTX}{MX} such that $\lambda\cdot rT=\gamma$. First we note that $\lambda=\inv{(rM)}\cdot R\gamma$, and therefore $\lambda$ is a natural transformation. Now, we observe that $\lambda\cdot(r\circ\eta)=\lambda\cdot rT\cdot\eta=\comp{\gamma}{\eta}=e$. Finally, since $\gamma$ and $rT$ are morphisms of monads, we have $\comp{rT}{\mu}=\comp{Rm}{(\ocomp{rT}{rT})}$ and $\comp{n}{(\ocomp{\gamma}{\gamma})}=\comp{\gamma}{\mu}$. By the middle-interchange law, we have $(\ocomp{\lambda}{\lambda})\cdot(\ocomp{rT}{rT})=\ocomp{(\comp{\lambda}{rT})}{(\comp{\lambda}{rT})}=\ocomp{\gamma}{\gamma}$. Consequently 
\begin{align*}
n\cdot(\ocomp{\lambda}{\lambda})\cdot(\ocomp{rT}{rT})&=\comp{n}{(\ocomp{\gamma}{\gamma})} \\ 
&=\gamma\cdot\mu\\
&=\comp{\lambda}{rT}\cdot\mu\\
&= \lambda\cdot\comp{Rm}{(\ocomp{rT}{rT})}.
\end{align*}
Note that $\ocomp{rT}{rT}=\comp{rTRT}{TrT}$. Since $rTRT$ and $TrT$ are epimorphisms, the result follows.
\end{proof}

The condition that $T$ should preserve epimorphisms may be compared with the notion of {\em compactification reflector} that is developed by Holgate in \cite{Hol09}. More specifically, a compactification reflector preserves the morphisms in a class $\mathcal{E}$ that is part of an ambient {\em factorisation system} $(\mathcal{E},\mathcal{M})$ on morphisms of \cat{C}. This feature is present in our examples on compactifications. We give the following observation.

\begin{prop}
\label{prop: separated monads are reflective within monads}
Let  $(R,r)$ be an epireflector on \cat{C} and let $MonE(\cat{C})$ be the category of those monads over \cat{C} that preserve epimorphisms and whose composition with $R$ gives a monad. Let $RMonE(\cat{C})$ be the subcategory of $MonE(\cat{C})$ consisting of those monads whose ranges are in $R(\cat{C})$. Then the embedding \map{I}{RMonE(\cat{C})}{MonE(\cat{C})} admits a left adjoint.
\end{prop}

\begin{proof}
We define the left adjoint $F$ as follows: $F(\trip{T})=\trip{R\circ T}$ and $F(\gamma)=R\gamma$ for any morphism \map{\gamma}{\trip{T}}{\trip{S}}. By Theorem \ref{thm: universality of the separation}, $F$ is well-defined and $Hom(F(\trip{T}),\trip{M})\cong Hom(\trip{T},I(\trip{M}))$.
\end{proof}

In the article \cite{Fak}, Fakir showed that under mild conditions, the category of idempotent monads is coreflective within the category of monads and monad morphisms. This construction was applied by Lambek and Rattray to describe - among other results, the Stone-\v{C}ech compactification and the Samuel compactification in \cite[Example 2 and Example 3]{LamRat}. Following Lambek and Rattray, Salbany used the construction to study what is known as {\em regular closure operators} (See \cite{Sal75} and \cite[Chapter 6]{DikTho}), a tool that is used to characterise epimorphisms in certain categories. (See \cite{DikKue} for instance.) The question that eventually arises is then that of comparing the reflective subcategories (as well as the various closure operators) associated to the idempotent cores of $T$ and $RT$. This question shall not be investigated here as it goes beyond the scope of the present paper.

\section{Fundamental examples}

In this section, we shall treat three separate and basic examples, and we shall give brief remarks about sobriety which will be helpful in clarifying certain aspects of the $T_0$ stable compactification. This will appear at the end of the first set of examples.

\subsection{Topological spaces and ultrafilter space monad}
We consider the ultrafilter space monad $\trip{U}=\mon{\ult}{\mu}{\eta}$ on \cat{Top}. For each topological space $X$,\ \map{\eta}{X}{\ult X} is a monomorphism (injective and continuous) which may fail to be an epimorphism (surjective and continuous), hence \ult\ is not idempotent. Given a surjective continuous map \map{f}{X}{Y}, each ultrafilter $\mathcal{G}\in\ult Y$ is the image of the ultrafilter with base $\{\inv{f}(G)\ |\ G\in\mathcal{G}\}$. Therefore \ult\ preserves epimorphisms. Finally, as already shown before, $\eta$ is a patch-dense and therefore is also a dense embedding.

\subsubsection{$T_0$ stable compactification}
We consider the $T_0$ epireflector $(R,r)$ on \cat{Top}. By Proposition \ref{prop: t0 quotient maps are proper} and Lemma \ref{lem: lifting a reflector as a monad}, the composition $R\ult$ forms a monad on \cat{Top}. As shown in Section \ref{The ultrafilter space monad}, $R\ult\cong\Sigma$. In fact, there is a unique monad morphism \map{\lambda}{\trip{R\circ U}}{\trip{S}} such that $\comp{\lambda}{r\ult}=\alpha$, where $\alpha$ is as defined in Proposition \ref{prop: monad morphism from prime to open prime}. With BUT, it can be shown as in Proposition \ref{prop: reflection of prime monad is open prime monad} that $\lambda$ is an homeomorphism. We make the following additional observations. (See \cite{Sal00,Smy,Sim}.)

\begin{prop}
\label{prop: salbany compactification}
For each continuous map \map{f}{X}{Y} where $(Y,\tau,\pi)$ is an \trip{U}-algebra, there is a unique map \map{\varphi}{\ult X}{Y} in \powcat{\cat{Top}}{U} such that $\comp{\varphi}{\eta}=f$.
\end{prop}

\begin{prop}
\label{prop: stable compactification}
For each continuous map \map{f}{X}{Y} where $(Y,\tau)$ is stably compact, there is a unique proper map \map{\varphi}{R\ult X}{Y} such that $\comp{\varphi}{(r\ult\cdot\eta)}=f$. The unit $r\ult\cdot\eta$ is a dense embedding if and only if $X$ is a $T_0$ space.
\end{prop}

\begin{prop}
\label{prop: stably compact space fails to be a reflective inside Top_0}
$R\eta$ is not an epimorphism in $\cat{Top}_0$.
\end{prop}

\begin{proof}
Since the inclusion $\cat{Stb}\to\cat{Top}_0$ is not full, \cat{Stb} is not reflective\footnote{At the end of the subsection, a different proof is given.} in $\cat{Top}_0$.
\end{proof}

\subsubsection{Prime closed filter compactification}

This method of compactification was the object of a thorough investigation by Bentley and Herrlich in the article \cite{BenHer}. The reader can consult \cite{Ban81}, \cite{Bru79}, \cite{Sal84} and \cite{Wyl84}, as well as the bibliography of \cite{BenHer} for an exhaustive reference. As we shall see here, the resulting construction is essentially that of Simmons' prime open filter monad. The monadic aspect of this construction was not explored in \cite{BenHer} and this is manifested in the absence of any attempt to define the multiplication.

Following \cite{BenHer}, let $\B_X$ denote the collection of all closed subsets of a space $X$. A {\em closed filter} on $X$ is a proper filter $\F$ on $X$ such that $\F\subseteq\B_X$, and \F\ is said to be a {\em prime closed filter} if \F\ is a closed filter and for any union of closed subsets $H\cup G\in\F$, either $H\in\F$ or $G\in\F$. The set of all prime closed filters on $X$ shall be denoted by $\prim X$ and \prim\ induces a monad $\trip{P}=(\prim,n,d)$ on \cat{Top} in the following manner: a base for closed sets on $\prim X$ is the collection \set{H^{\times}\ |\ H\in\B_X}, with $H^{\times}=\set{\F\in\prim X\ |\ H\in\F}$. The natural transformations $n$ and $d$ are respectively defined by
\begin{center}
$n_X(\mathfrak{X})=\{H\in\B_X\setminus\{\emptyset\}\ |\ H^{\times}\in\mathfrak{X}\}$ and $d_X(x)=\{C\in\B_X\ |\ x\in C\}$.
\end{center}
For each continuous map \map{f}{X}{Y}, we have $\prim f(\F)=\set{B\in\B_Y\ |\ \inv{f}(B)\in\F}$.

As shown in \cite[Proposition 3]{BenHer} we have a surjective map \map{\alpha}{\ult X}{\prim X} given by $\F\mapsto\F\cap\B_X$. (BUT is assumed.) This map is continuous for each space $X$ and is a monad morphism (Proposition \ref{prop: monad morphism from prime to open prime}). Now, since $\prim X$ is $T_0$ (\cite[Remarks 2.1]{BenHer}), there is a unique morphism of monads \map{\lambda}{\trip{R\circ U}}{\trip{P}} such that $\comp{\lambda}{r\ult}=\alpha$. (Here $\trip{R}=(R,r)$ is still the $T_0$ epireflector.) As in the previous case, $\lambda$ becomes an homeomorphism. This can be made more transparent as follows.

\begin{lem}
\label{lem: T_0 reflection with closed sets}
If $(R,r)$ is the $T_0$ epireflector and \map{r}{X}{RX} a $T_0$ quotient map, then for any pair of points $x$ and $y$, we have $r(x)=r(y)$ if and only if for any closed set $C$, $x\in C$ whenever $y\in C$ and vice versa. In particular $r\ult(\F)=r\ult(\G)$ if and only if $\F\cap\B_X=\G\cap\B_X$.
\end{lem}

\begin{prop}
\label{prop: isomorphism of monads}
Granted that BUT holds, there is a monad isomorphism between the prime open filter monad $\trip{S}$ and the prime closed filter monad \trip{P}.
\end{prop}

\subsubsection{\v{C}ech-Stone compactification}
In \cat{Top}, we consider the full subcategories \cat{Haus} of Hausdorff spaces and \cat{CHaus} of compact Hausdorff spaces. The epireflector onto \cat{Haus} shall be denoted by $(H,h)$. 

\begin{lem}
\label{lem: Reta is dense in haus}
$H\eta$ is an epimorphism in \cat{Haus}.
\end{lem}

\begin{proof}
In the diagram $\comp{h\ult}{\eta}=\comp{H\eta}{h}$, $h\ult$ and $h$ are surjective and $\eta$ is dense. Consequently $H\eta$ is dense, i.e. it is an epimorphism in \cat{Haus}.
\end{proof}

\begin{lem}
\label{lem: Hausdorff reflection of stably compact}
The reflection $H\ult X$ of $\ult X$ is compact and Hausdorff.
\end{lem}
Since any compact and Hausdorff space is stably compact, the space $H\ult X$ is then a \trip{U}-algebra. It follows from Theorem \ref{thm: preservation of free algebras is enough}, Corollary \ref{cor: strong lifting of a reflector as a monad} and Proposition \ref{prop: properties of the rt-algebras as reflective subcategories} that

\begin{prop}
\label{prop: Czech-Stone compactification}
$H\ult$ is a reflector that is left adjoint to the full inclusion $\cat{CHaus}\to\cat{Top}$. Furthermore \cat{CHaus} is reflective in \cat{Haus} and \powcat{\cat{Top}}{U}.
\end{prop}
Now, consider the prime open filter monad \trip{S}. It is straightforward to check that $He$ is an epimorphism in $\cat{Haus}$ and that the reflection $H\Sigma X$ of $\Sigma X$ is compact and Hausdorff. Consequently

\begin{prop}
\label{prop: CHaus reflective in Stb}
\cat{CHaus} is reflective in \cat{Stb}.
\end{prop}
The well-known fact that any continuous map in \cat{CHaus} must be proper is represented through the factorisation $\cat{CHaus}\to\cat{Stb}\to\cat{Top}$. With Lemma \ref{lem: CHaus coreflective in Stb}, we conclude that \cat{CHaus} is simultaneously reflective and coreflective in \cat{Stb}. (The localic version of this is discussed by Escard\'o in \cite{Esc}.) The fact that \cat{CHaus} is reflective in \cat{Stb} is shown in \cite[Section 5.6]{Low} in the context of lax algebras. That $H\ult\cong H\Sigma$ is an application of Theorem \ref{thm: universality of the separation}. Thus $H\ult$ represents the \v{C}ech-Stone compactification $\beta$. In the absence of BUT, it can be concluded that there are morphisms of monads that arise from different separation of \ult, making the following diagram commutative.
$$\xymatrix{
& & \Sigma \ar[rd] &\\
\ult \ar[r] & R\ult \ar[dr]\ar[ur] & & \beta \\
& & \prim \ar[ru] &
}$$

\paragraph*{Remarks about sobriety.}
Let \map{(\sob,s)}{\cat{Top}}{\cat{Top}} be the sobrification reflector. Since $\ult X$ is weakly sober for each space $X$, one has $\sob\ult X\cong R\ult X$ where $(R,r)$ is the $T_0$ reflector. Thus \sob\ does not improve $\ult$ in a way that is different from $R$. The functor \sob\ may be better understood as a completion than a separation condition, especially in the context of $T_0$ spaces. (See \cite{KueBru} and \cite{RHof76}.) The embedding \map{s}{X}{\sob X} is $b$-dense hence an epimorphism in $\cat{Top}_0$. Note that since $\sob X$ is naturally isomorphic to the set of completely prime open filters on $X$ with an appropriate topology (\cite{Joh}), $\Sigma X\ncong\sob X$. This also shows Proposition \ref{prop: stably compact space fails to be a reflective inside Top_0}: that \map{e}{X}{\Sigma X} is not $b$-dense hence not an epimorphism in $\cat{Top}_0$, for otherwise $\Sigma X$ would act as a sobrification of $X$. (See explanation in \cite[Proposition 3.1.2]{RHof76}.)

\subsection{Quasi-uniform spaces and completion monad}
We denote by \cat{QUnif} the category of quasi-uniform spaces and quasi-uniformly continuous maps. Following Salbany (\cite{Sal82}), we consider the set $\cau X$ of all Cauchy filters on a quasi-uniform space $(X,\qu)$ and endow it with the quasi-uniformity $\qu^*$ having as a basis the set \set{U^*\ |\ U\in\qu} where $(\F,\G)\in U^*$ if and only if there are $(A,B)\in\F\times\G$ such that $A\times B\subseteq U$. For each Cauchy filter \F,\ $\cau f(\F)$ has as a basis the set $\set{G\ |\ \inv{f}(G)\in\F}$. Thus \cau\ becomes an endofunctor on \cat{QUnif} that induces a monad $\trip{C}=(\cau,\mu,\eta)$ with
\begin{center}
$\mu_X(\mathfrak{X})=\set{G\ |\ G^*\in\mathfrak{X}}$ and $\eta_X(x)=\set{A\ |\ x\in A}$,
\end{center}
and where $x\in X$, $\mathfrak{X}\in\cau\cau X$ and $G^*=\set{\F\ |\ G\in\F}$. Apart from \cite{Sal82}, this monad has also been discussed by Wyler in \cite{Wyl74}. We consider the ``separation'' epireflector \map{(R,r)}{\cat{QUnif}}{\cat{QUNif}} with $R(\cat{QUnif})=\cat{QUnif}_0$. (See \cite{Sal82}.) We gather from \cite{Sal82} that $\eta(X)$ is dense with respect to the topology induced by $\qu^*\vee(\qu^*)^{-1}$ and that the algebras \powcat{\cat{QUnif}}{C} are the complete quasi-uniform spaces. 

\begin{lem}
\label{lem: completion monad not idempotent}
The monad \trip{C} is not idempotent.
\end{lem}

\begin{proof}
Note first that epimorphisms in \cat{QUnif} are precisely the surjective quasi-uniformly continuous maps (\cite{DikKue}). The unit $\eta$, which is a monomorphism, fails to be surjective in general.
\end{proof}

\begin{lem}
\label{lem: Reta epi in Qunif0}
$R\eta$ is an epimorphism in $\cat{QUnif}_0$ and $R\cau (X,\qu)$ is a separated complete quasi-uniform space for each quasi-uniform space $(X,\qu)$.
\end{lem}

\begin{proof}
As indicated previously (\cite{Sal82}), $\eta$ is dense with respect to $\qu^*\vee(\qu^*)^{-1}$. Therefore the compositions in the diagram $\comp{r\cau}{\eta}=\comp{R\eta}{r}$ are dense. This shows that $R\eta$ is dense with respect to $(\qu^*)^{s}\vee\inv{((\qu^*)^{s})}$, hence an epimorphism in $\cat{QUnif}_0$.
\end{proof}
Consequently $R\cau$ is a reflector and the complete separated quasi-uniform spaces form a reflective subcategory of $\cat{QUnif}$, $\cat{QUnif}_0$ and \powcat{\cat{QUnif}}{C}. These facts which are indicated or shown in \cite{Sal82} have substantially motivated the study of separation with respect to a monad. The separated completion monad $R\cau$ is also investigated in \cite{CleHof} where quasi-uniform spaces are viewed as lax proalgebras. Various epimorphisms were studied by Dikranjan and K\"unzi in relation to the different conditions of separation in \cat{QUnif}. These various separation conditions could not be considered here due to the author's limited expertise.

\subsection{Frames and ideal functor}
We refer the reader to \cite{Joh}, \cite{Ban81} and \cite{Ban90,BanBru} for general background on frames and the ideal functor. A complete lattice $L$ is a frame if $a\wedge(\bigvee_Ib_i)=\bigvee_I(a\wedge b_i)$ for any $a\in L$ and $\set{b_i\ |\ i\in I}\subseteq L$. A map \map{f}{L}{M} is a frame homomorphism if it preserves arbitrary joins and finite meets. We denote by \cat{Frm} the category of frames and frame homomorphisms. The symbol \id\ will denote the ideal functor that takes each frame $L$ to its set of ideals $\id L$ which is a frame, and takes each frame homomorphism \map{f}{L}{M} to a frame homomorphism $\map{\id f}{\id L}{\id M}: I\mapsto \bigcup \downarrow\{f(a)\ |\ a\in I\}=\{b\leq f(a)\ |\ a\in I\}$. The functor \id\ is part of a comonad $\trip{I}=\mon{\id}{c}{\sigma}$ on \cat{Frm} where $\sigma_L:\id L\to L: I\mapsto\bigvee I$ for all $I\in\id L$, and $\map{c_L}{\id L}{\id\id L}:J\mapsto\{I\in\id L\ |\ \sigma_L(I)\in J\}$.

The way below relation $\ll$ is defined as for stably compact spaces: $a\ll b$ if any arbitrary subset whose join is above $b$ admits a finite subset whose join is above $a$. A frame $L$ is stably continuous if the way below relation $\ll$ is continuous (\cite{BanBru,Gie,Joh}), i.e. $a=\bigvee\set{b\ |\ b\ll a}$ for all $a\in L$, and finitely multiplicative. (Hence $L$ is compact.) For a comprehensive analysis of continuous lattices in general, we invite the reader to consult \cite{Gie}. A frame homomorphism is proper if it preserves the way below relation. The category of stably continuous frames together with proper frame maps is denoted by \cat{SCont}.

Although it is shown in \cite{BanBru} that \trip{I} is induced from stably continuous frames in a sense that $\powcat{\cat{Frm}}{I}\cong\cat{SCont}$, the co-unit $\sigma$ and co-multiplication $c$ can be described independently as above by only using the ambient category \cat{Frm}. Our objective is to give an alternative description of the regular - and consequently completely regular, compactification of frames.

A frame $L$ is said to be regular if for all $a\in L$, $a=\bigvee\{c\ |\ c\prec a\}$, where $a\prec b$\ \ if\ \ $b\vee a^*=1$ with $a^*=\bigvee\{x\in L\ |\ x\wedge a=0\}$. If $Reg L$ is the largest regular subframe of $L$ and \map{r_L}{Reg L}{L} the inclusion frame homomorphism, then \map{(Reg,r)}{\cat{Frm}}{\cat{Frm}} is a monocoreflector.

\begin{lem}
\label{lem: compact regular are stably compact}
For each frame $L$,\ the coreflection $Reg\id L$ is compact and regular, hence stably continuous. $Reg\sigma$ is a monomorphism in the category of regular frames $Reg(\cat{Frm})=\cat{RegFrm}$.
\end{lem}

\begin{proof}
The first statement is trivial. Monomorphisms in \cat{RegFrm} are precisely the dense frame homomorphisms, i.e. those frame homomorphisms that reflect $0$. Let us first note that if $\sigma_L(I)=0$ where $I\in\id L$, then $\bigvee I=0$ and so $I=\set{0}$. Therefore $\sigma$ is dense. As in the previous cases (Lemma \ref{lem: Reta is dense in haus} and Lemma \ref{lem: Reta epi in Qunif0}), since the compositions $\comp{r_L}{Reg\sigma}=\comp{\sigma_L}{r\id}$ are dense and $r_L$ is a monomorphism, $Reg\sigma$ is dense.
\end{proof}
It then follows that $Reg\id$ is a coreflector (idempotent comonad), i.e. the category of compact regular frames and frame homomorphisms \cat{KRegFrm} is coreflective in \cat{Frm}, \cat{SCont} and \cat{RegFrm}. The fact that any frame homomorphism between compact regular frames is proper is syntactically encoded through the factorisation $\cat{KRegFrm}\to\cat{Stb}\to\cat{Frm}$. Here also \cat{KRegFrm} is simultaneously reflective and coreflective in \cat{SCont}. (See \cite{Esc}.) We end with a final observation.

\begin{prop}
\label{prop: ideal functor is a compactification}
The functor \id\ preserves monomorphisms. Consequently if $(R,r)$ is a monocoreflector on \cat{Frm} such that $R\id$ is a comonad, then any comonad morphism \map{\gamma}{T}{\id} with $T(\cat{Frm})\subseteq R(\cat{Frm})$ factors uniquely through a comonad morphism \map{\lambda}{T}{R\id}.
\end{prop}

\paragraph*{Acknowledgements} I thank Dirk Hofmann for insightful discussions on the ultrafilter space monad, and for directing me to the work of H. Simmons in \cite{Sim} and the article \cite{Bez} by Bezhanishvili et al. I am grateful to Themba Dube and Anneliese Schauerte for having contributed to my knowledge on the ideal functor and on stably compact and regular frames. Finally, I am grateful to the anonymous referees whose generous suggestions have contributed to a better exposition of the paper.



\end{document}